\documentclass{gtpart}  
\usepackage{pinlabel}
\usepackage[all]{xy}
\title{Geometric embedding properties \\ of Bestvina-Brady subgroups}

\author{Hung Cong Tran} 
\givenname{Hung Cong}
\surname{Tran}
\email{hung.tran@uga.edu}
\address{Department of Mathematics\\
 The University of Georgia\\
1023 D. W. Brooks Drive\\
Athens, GA 30605\\
United States}

\arxivreference{arXiv:1606.00539}  %%% please supply if 
\arxivpassword{k6wgm}   %%% paper is in the arXiv

\newtheorem{thm}{Theorem}[section]
\newtheorem{prop}[thm]{Proposition}
\newtheorem{lem}[thm]{Lemma}
\newtheorem{cor}[thm]{Corollary}

\theoremstyle{definition}
\newtheorem{defn}[thm]{Definition}
\newtheorem{rem}[thm]{Remark}

\newcommand{\abs}[1]{\lvert{#1}\rvert}

\renewcommand{\bar}[1]{\overline{#1}}

\newcommand{\bigset}[2]{ \bigl\{ \, {#1} \bigm| {#2} \, \bigr\} }

\DeclareMathOperator{\CAT}{CAT}

\DeclareMathOperator{\Div}{Div}

\begin{document}

%\title[Geometric embedding properties of Bestvina-Brady subgroups]{Geometric embedding properties of Bestvina-Brady subgroups}

%\author{Hung Cong Tran}
%\address{Department of Mathematics\\
 %The University of Georgia\\
%1023 D. W. Brooks Drive\\
%Athens, GA 30605\\
%USA}
%\email{hung.tran@uga.edu}

%\date{\today}

\begin{abstract}
We compute the relative divergence of right-angled Artin groups with respect to their Bestvina-Brady subgroups and the subgroup distortion of Bestvina-Brady subgroups. We also show that for each integer $n\geq 3$, there is a free subgroup of rank $n$ of some right-angled Artin group whose inclusion is not a quasi-isometric embedding. The corollary answers the question of Carr about the minimum rank $n$ such that some right-angled Artin group has a free subgroup of rank $n$ whose inclusion is not a quasi-isometric embedding. It is well-known that a right-angled Artin group $A_\Gamma$ is the fundamental group of a graph manifold whenever the defining graph $\Gamma$ is a tree with at least 3 vertices. We show that the Bestvina-Brady subgroup $H_\Gamma$ in this case is a horizontal surface subgroup.
\end{abstract}

%\subjclass[2000]{%
%20F67, % Hyperbolic groups and nonpositively curved groups
%20F65} % Geometric group theory
\maketitle

\section{Introduction}
For each $\Gamma$ a finite simplicial graph the associated \emph{right-angled Artin group} $A_{\Gamma}$ has generating set $S$ the vertices of $\Gamma$, and relations $st = ts$ whenever $s$ and $t$ are adjacent vertices. If $\Gamma$ is non-empty, there is a homomorphism from $A_\Gamma$ onto the
integers, that takes every generator to $1$. The \emph{Bestvina-Brady subgroup} $H_\Gamma$ is defined to be the kernel of this homomorphism. 

Bestvina-Brady subgroups were introduced by Bestvina-Brady in \cite{MR1465330} to study the finiteness properties of subgroups of right-angled Artin groups. One result in \cite{MR1465330} is that the Bestvina-Brady subgroup $H_{\Gamma}$ is finitely generated iff the graph $\Gamma$ is connected. This fact is a motivation to study the geometric connection between a right-angled Artin group and its Bestvina-Brady subgroup. More precisely, we examine the relative divergence of right-angled Artin groups with respect to their Bestvina-Brady subgroups and the subgroup distortion of Bestvina-Brady subgroups (see the following theorem).

\begin{thm}
\label{main}
Let $\Gamma$ be a connected, finite, simplicial graph with at least 2 vertices. Let $A_{\Gamma}$ be the associated right-angled Artin group and $H_{\Gamma}$ the Bestvina-Brady subgroup. Then the relative divergence $Div(A_{\Gamma}, H_{\Gamma})$ and the subgroup distortion $Dist^{H_{\Gamma}}_{A_{\Gamma}}$ are both linear if $\Gamma$ is a join graph. Otherwise, the relative divergence $Div(A_{\Gamma}, H_{\Gamma})$ and the subgroup distortion $Dist^{H_{\Gamma}}_{A_{\Gamma}}$ are both quadratic. 
\end{thm}

In the above theorem, we can see that the relative divergence $Div(A_{\Gamma}, H_{\Gamma})$ and the subgroup distortion $Dist^{H_{\Gamma}}_{A_{\Gamma}}$ are equivalent. In general, we showed that the relative divergence is always dominated by the subgroup distortion for any pair of finitely generated groups $(G, H)$, where $H$ is a normal subgroup of $G$ such that the quotient group $G/H$ is an infinite cyclic group (see
Proposition \ref{p1}).

Carr \cite{MC} proved that non-abelian two-generator subgroups of right-angled Artin groups are quasi-isometrically embedded free groups. In his paper, he also showed an example of a distorted free subgroup of a right-angled Artin group. However, the minimum rank $n$ such that some right-angled Artin group has a free subgroup of rank $n$ whose inclusion is not a quasi-isometric embedding was still unknown (see \cite{MC}). A corollary of Theorem \ref{main} answered this question (see the following corollary).

\begin{cor}
\label{hch}
For each integer $n\geq 3$, there is a right-angled Artin group containing a free subgroup of rank $n$ whose inclusion is not a quasi-isometric embedding.
\end{cor} 

We remark that a special case of Theorem \ref{main} can also be derived as a consequence of previous work by Hruska--Nguyen (\cite{H-N}) on distortion of surfaces in graph manifolds. Hruska--Nguyen showed that every virtually embedded horizontal surface in a $3$--dimensional graph manifold has quadratic distortion. After learning about this result, the we proved the following theorem, which implies that many Bestvina--Brady subgroups are also horizontal surface subgroups.

\begin{thm}
\label{hzs}
If $\Gamma$ is a finite tree with at least 3 vertices, then the associated right-angled Artin group $A_\Gamma$ is a fundamental group of a graph manifold and the Bestvina-Brady subgroup $H_\Gamma$ is a horizontal surface subgroup.
\end{thm} 

It is well-known that a right-angled Artin group $A_\Gamma$ is the fundamental group of a graph manifold whenever the defining graph $\Gamma$ is a tree with at least 3 vertices. However, the fact that the Bestvina-Brady subgroup $H_\Gamma$ is a horizontal subgroup does not seem to be recorded in the literature. With the use of Theorem \ref{hzs}, we see that Theorem \ref{main} can be viewed as a generalization of a special case of the quadratic distortion theorem of Hruska--Nguyen. Moreover, Theorem \ref{hzs} combined with the Hruska--Nguyen theorem gives an alternative proof of Corollary \ref{hch}.

\subsection*{Acknowledgments}
I would like to thank Prof.~Christopher Hruska and Hoang Nguyen for very helpful comments and suggestions, especially their help on the proof of Theorem \ref{hzs}. Hoang Nguyen helped me a lot with the background of graph manifolds and horizontal surfaces. He also showed me his joint work with Hruska on subgroup distortion of horizontal surface subgroups and their work is a motivation for this article. I also thank the referee for advice that improved the exposition of the paper.

\section{Right-angled Artin groups and Bestvina-Brady subgroups}

\begin{defn}
Given a finite simplicial graph $\Gamma$, the associated \emph{right-angled Artin group} $A_{\Gamma}$ has generating set $S$ the vertices of $\Gamma$, and relations $st = ts$ whenever $s$ and $t$ are adjacent vertices.

Let $S_1$ be a subset of $S$. The subgroup of $A_{\Gamma}$ generated by $S_1$ is a right-angled Artin group $A_{\Gamma_1}$, where $\Gamma_1$ is the induced subgraph of $\Gamma$ with vertex set $S_1$ (i.e.~$\Gamma_1$ is the union of all edges of $\Gamma$ with both endpoints in $S_1$). The subgroup $A_{\Gamma_1}$ is called a \emph{special subgroup} of $A_{\Gamma}$. 
\end{defn}
\
\begin{defn}
Let $\Gamma$ be a finite simplicial graph with the set $S$ of vertices. Let $T$ be a torus of dimension $\abs{S}$ with edges labeled by the elements of $S$. Let $X_{\Gamma}$ denote the subcomplex of $T$ consisting of all faces whose edge labels span a complete subgraph in $\Gamma$ (or equivalently, mutually commute in $A_{\Gamma}$). $X_{\Gamma}$ is called the \emph{Salvetti complex}.
\end{defn}

\begin{rem}
The fundamental group of $X_{\Gamma}$ is $A_{\Gamma}$. The universal cover $\tilde{X}_{\Gamma}$ of $X_{\Gamma}$ is a $\CAT(0)$ cube complex with a free, cocompact action of $A_{\Gamma}$. Obviously, the 1--skeleton of $\tilde{X}_{\Gamma}$ is the Cayley graph of $A_{\Gamma}$ with respect the generating set $S$.
\end{rem}

\begin{defn}
Let $\Gamma$ be a finite simplicial graph. Let $\Phi\!: A_{\Gamma} \to \Z$ be an epimorphism which sends all the generators of $A_{\Gamma}$ to $1$ in $\Z$. The kernel $H_{\Gamma}$ of $\Phi$ is called the \emph{Bestvina-Brady subgroup}.
\end{defn}

\begin{rem}
\label{r1}
There is a natural continuous map $f\!: X_{\Gamma} \to S^{1}$ which induces the homomorphism $\Phi\!: A_{\Gamma} \to \Z$. Moreover, it is not hard to see that the lifting map $\tilde{f}\!: \tilde{X}_{\Gamma} \to \R$ is an extension of $\Phi$.
\end{rem}

\begin{thm}[Bestvina-Brady \cite{MR1465330} and Dicks-Leary \cite{MR1605948}]
\label{th1}
Let $\Gamma$ be a finite simplicial graph. The Bestvina-Brady subgroup $H_{\Gamma}$ is finitely generated iff $\Gamma$ is connected. Moreover, the set $T$ of all elements of the form $st^{-1}$ whenever $s$ and $t$ are adjacent vertices form a finite generating set for $H_{\Gamma}$. Moreover, if $\Gamma$ is a tree with $n$ edges, then the Bestvina-Brady subgroup $H_{\Gamma}$ is a free group of rank $n$.
\end{thm}

\begin{defn}
Let $\Gamma_1$ and $\Gamma_2$ be two graphs, the \emph{join} of $\Gamma_1$ and $\Gamma_2$ is a graph obtained by connecting every vertex of $\Gamma_1$ to every vertex of $\Gamma_2$ by an edge.

Let $J$ be a complete subgraph of $\Gamma$ which decomposes as a nontrivial join. We call $A_J$ a \emph{join subgroup} of $A_{\Gamma}$.
\end{defn}

Let $\Gamma$ be a finite simplicial graph with the vertex set $S$ and let $g$ an element of $A_\Gamma$. A \emph{reduced word} for $g$ is a minimal length word in the free group $F(S)$ representing $g$. Given an arbitrary word representing $g$, one can obtain a reduced word by a process of ``shuffling'' (i.e. interchanging commuting elements) and canceling inverse pairs. Any two reduced words for $g$ differ only by shuffling. 
For an element $g \in A_\Gamma$, a cyclic reduction of $g$ is a minimal length element of the conjugacy class of $g$. If $w$ is a reduced word representing $g$, then we can find a cyclic reduction $\bar{g}$ by shuffling commuting generators in $w$ to get a maximal length word $u$ such that $w=u\bar{w}u^{-1}$. In particular, $g$ itself is \emph{cyclically reduced} if and only if every shuffle of $w$ is cyclically reduced as a word in the free group $F(S)$.

\section{Relative divergence, geodesic divergence, and subgroup distortion}
Before we define the concepts of relative divergence, geodesic divergence, and subgroup distortion, we need to build the notions of domination and equivalence. These notions are the tools to measure the relative divergence, geodesic divergence, and subgroup distortion.
\begin{defn}
Let $\mathcal{M}$ be the collection of all functions from $[0,\infty)$ to $[0,\infty]$. Let $f$ and $g$ be arbitrary elements of $\mathcal{M}$. \emph{The function $f$ is dominated by the function $g$}, denoted \emph{$f\preceq g$}, if there are positive constants $A$, $B$, $C$ and $D$ such that $f(x)\leq Ag(Bx)+Cx$ for all $x>D$. Two function $f$ and $g$ are \emph{equivalent}, denoted \emph{$f\sim g$}, if $f\preceq g$ and $g\preceq f$. %\emph{The function $f$ is strictly dominated by the function $g$}, denoted \emph{$f\prec g$}, if $f$ is dominated by $g$ and they are not equivalent.

\end{defn}

\begin{rem}
%The relations $\preceq$ and $\prec$ are transitive. The relation $\sim$ is an equivalence relation on the set $\mathcal{M}$. 
%Let $f$ and $g$ be two polynomial functions in the family $\mathcal{M}$. We observe that $f$ is dominated by $g$ iff the degree of $f$ is less than or equal to the degree of $g$ and they are equivalent iff they have the same degree. All exponential functions of the form $a^{bx+c}$, where $a>1, b>0$ are equivalent. Therefore, 
A function $f$ in $\mathcal{M}$ is \emph{linear, quadratic or exponential...} if $f$ is respectively equivalent to any polynomial with degree one, two or any function of the form $a^{bx+c}$, where $a>1, b>0$.
\end{rem}

\begin{defn}
Let $\{\delta^n_{\rho}\}$ and $\{\delta'^n_{\rho}\}$ be two families of functions of $\mathcal{M}$, indexed over $\rho \in (0,1]$ and positive integers $n\geq 2$. \emph{The family $\{\delta^n_{\rho}\}$ is dominated by the family $\{\delta'^n_{\rho}\}$}, denoted \emph{$\{\delta^n_{\rho}\}\preceq \{\delta'^n_{\rho}\}$}, if there exists constant $L\in (0,1]$ and a positive integer $M$ such that $\delta^n_{L\rho}\preceq \delta'^{Mn}_{\rho}$. Two families $\{\delta^n_{\rho}\}$ and $\{\delta'^n_{\rho}\}$ are \emph{equivalent}, denoted \emph{$\{\delta^n_{\rho}\}\sim \{\delta'^n_{\rho}\}$}, if $\{\delta^n_{\rho}\}\preceq \{\delta'^n_{\rho}\}$ and $\{\delta'^n_{\rho}\}\preceq \{\delta^n_{\rho}\}$.
\end{defn}

\begin{rem}
%The relations $\preceq$ and $\prec$ are transitive. The relation $\sim$ is an equivalence relation. 
%If $f$ is an element in $\mathcal{M}$, we could represent $f$ as a family $\{\delta^n_{\rho}\}$ for which $\delta^n_{\rho}=f$ for all $\rho$ and $n$. 
A family $\{\delta^n_{\rho}\}$ is dominated by (or dominates) a function $f$ in $\mathcal{M}$ if $\{\delta^n_{\rho}\}$ is dominated by (or dominates) the family $\{\delta'^n_{\rho}\}$ where $\delta'^n_{\rho}=f$ for all $\rho$ and $n$. The equivalence between a family $\{\delta^n_{\rho}\}$ and a function $f$ in $\mathcal{M}$ can be defined similarly. Thus, a family $\{\delta^n_{\rho}\}$ is linear, quadratic, exponential, etc if $\{\delta^n_{\rho}\}$ is equivalent to the function $f$ where $f$ is linear, quadratic, exponential, etc.
\end{rem}

\begin{defn}
Let $X$ be a geodesic space and $A$ a subspace of $X$. Let $r$ be any positive number.
\begin{enumerate}
\item $N_r(A)=\bigset{x \in X}{d_X(x, A)<r}$
\item $\partial N_r(A)=\bigset{x \in X}{d_X(x, A)=r}$ 
\item $C_r(A)=X-N_r(A)$.
\item Let $d_{r,A}$ be the induced length metric on the complement of the $r$--neighborhood of $A$ in $X$. If the subspace $A$ is clear from context, we can use the notation $d_r$ instead of using $d_{r,A}$. 
\end{enumerate}
\end{defn}

\begin{defn}
Let $(X,A)$ be a pair of metric spaces. For each $\rho \in (0,1]$ and positive integer $n\geq 2$, we define a function $\delta^n_{\rho}\!:[0, \infty)\to [0, \infty]$ as follows: 

For each $r$, let $\delta^n_{\rho}(r)=\sup d_{\rho r}(x_1,x_2)$ where the supremum is taken over all $x_1, x_2 \in \partial N_r(A)$ such that $d_r(x_1, x_2)<\infty$ and $d(x_1,x_2)\leq nr$. 
%\[\delta^n_{\rho}(r)=\sup\set{d_{C_{\rho r}(A)}(x_1, x_2)}{\text{$x_1$, $x_2$ $\in \partial N_r(A)$, $d_{C_r(A)}(x_1, x_2)< \infty$, and $d_S(x_1,x_2)\leq nr$}}.\]

%\[\delta^n_{\rho}(r)=\sup_{\substack{x_1,x_2\in \partial N_r(A)\\d_r(x_1, x_2)< \infty\\d_S(x_1,x_2)\leq nr}} d_{\rho r}(x_1, x_2).\]

The family of functions $\{\delta^n_{\rho}\}$ is \emph{the relative divergence} of $X$ with respect $A$, denoted $Div(X,A)$.
\end{defn}

We now define the concept of relative divergence of a finitely generated group with respect to a subgroup.
\begin{defn} 
Let $G$ be a finitely generated group and $H$ its subgroup. We define \emph{the relative divergence} of $G$ with respect to $H$, denoted \emph{$Div(G,H)$}, to be the relative divergence of the Cayley graph $\Gamma(G,S)$ with respect to $H$ for some finite generating set $S$. 
\end{defn}

\begin{rem}
The concept of relative divergence was introduced by the author in \cite{MR3361149} with the name upper relative divergence. The relative divergence of geodesic spaces is a pair quasi-isometry invariant concept. This implies that the relative divergence on a finitely generated group does not depend on the choice of finite generating sets. 
\end{rem}

\begin{defn}
The \emph{divergence} of a bi-infinite geodesic $\alpha$, denoted $\Div_{\alpha}$, is a function $g\!: (0,\infty)\to(0,\infty)$ which for each positive number $r$ the value $g(r)$ is the infimum on the lengths of all paths outside the open ball with radius $r$ about $\alpha(0)$ connecting $\alpha(-r)$ and $\alpha(r)$.
\end{defn}

The following lemma is deduced from the proof of Corollary 4.8 in \cite{MR2874959}.

\begin{lem}
\label{l4}
Let $\Gamma$ be a connected, finite, simplicial graph with at least 2 vertices. Assume that $\Gamma$ is not a join. Let $g$ be a cyclically reduced element in $A_{\Gamma}$ that does not lie in any join subgroup. Then the divergence of bi-infinite geodesic $\cdots ggggg\cdots$ is at least quadratic. 
\end{lem}

\begin{defn}
Let $G$ be a group with a finite generating set $S$ and $H$ a subgroup of $G$ with a finite generating set $T$. 
The \emph{subgroup distortion} of $H$ in $G$ is the function $Dist^H_G\!:(0,\infty)\to(0,\infty)$ defined as follows:
\[Dist^H_G(r)=\max \bigset{\abs{h}_T}{h\in H, \abs{h}_S\leq r}.\] 
\end{defn}

\begin{rem}
It is well-known that the concept of distortion does not depend on the choice of finite generating sets.
\end{rem}
\section{Connection between subgroup distortion and relative divergence} 

\begin{lem}
\label{l3}
Let $H$ be a finitely generated group with finite generating set $T$ and $\phi$ in $Aut(H)$. Let $G=\langle H,t/tht^{-1}=\phi(h)\rangle$ and $S=T\cup\{t\}$. Then:
\begin{enumerate}
\item All element in $G$ can be written uniquely in the from $ht^n$ where $h$ is a group element in $H$.
\item The set $S$ is a finite generating set of $G$ and $d_S(ht^m,h't^n)\geq \abs{m-n}$, and $d_S(ht^m,Ht^n)=\abs{m-n}$. 
\end{enumerate} 
\end{lem}
\begin{proof}
The statement (1) is a well-known and we only need to prove statement (2). Let $\psi$ be the map from $G$ to $\Z$ by sending element $t$ to 1 and each generator in $T$ to 0. It is not hard to see that $\psi$ is a group homomorphism. We first show that the absolute value of $\psi(g)$ is at most the length of $g$ with respect to $S$ for each group element $g$ in $G$. In fact, let $w_1t^{n_1}w_2t^{n_2}\cdots w_kt^{n_k}$ be the shortest word in $S$ that represents $g$, where each $w_i$ is a word in $T$. Therefore,
 \[\psi(g)=n_1+n_2+\cdots+n_k\]
and 
\[\abs{g}_S= \bigl(\ell(w_1)+\ell(w_2)+\cdots+\ell(w_k)\bigr)+\bigl(\abs{n_1}+\abs{n_2}+\cdots+\abs{n_k}\bigr).\]
This implies that the absolute value of $\psi(g)$ is at most the length of $g$ with respect to $S$. The distance between two elements $ht^m$ and $h't^n$ is the length of the group element $g=(ht^m)^{-1}h't^n$. Obviously, $\psi(g)=n-m$. Therefore, the distance between two elements $ht^m$ and $h't^n$ is at least $\abs{m-n}$. This fact directly implies that the distance between $ht^m$ and any element in $Ht^n$ is at least $\abs{m-n}$. Also, $ht^n$ is an element in $Ht^n$ and the distance between $ht^m$, $ht^n$ is at most $\abs{m-n}$. Therefore, the distance between $ht^m$ and $Ht^n$ is at exactly $\abs{m-n}$. 
\end{proof}

\begin{lem}
\label{l5}
Let $H$ be a finitely generated group with finite generating set $T$ and $\phi$ in $Aut(H)$. Let $G=\langle H,t/tht^{-1}=\phi(h)\rangle$ and $S=T\cup\{t\}$. Let $n$ be an arbitrary positive integer and $x$, $y$ be two points in $\partial N_n(H)$. Then there is path outside $N_n(H)$ connecting $x$ and $y$ iff the pair $(x,y)$ is either of the form $(h_1t^n, h_2t^n)$ or $(h_1t^{-n}, h_2t^{-n})$ where $h_1$ and $h_2$ are elements in $H$. 
\end{lem}
\begin{proof}
By Lemma \ref{l3}, the pair $(x,y)$ must be of the form $(h_1t^{m_1}, h_2t^{m_2})$ where $\abs{m_1}=\abs{m_2}=n$. We first assume that $m_1m_2<0$. Let $\gamma$ be an arbitrary path connecting $x$ and $y$. By Lemma \ref{l3}, we observe that if two vertices $ht^m$ and $h't^{m'}$ of $\gamma$ are consecutive, then $\abs{m-m'}\leq 1$. Therefore, there exists a vertex of $\gamma$ that belongs to $H$. Thus, there is no path outside $N_n(H)$ connecting $x$ and $y$. 

If $m_1=m_2$, then $x$ and $y$ both lie in the same coset $t_{m_1}H$. Therefore, there is a path $\alpha$ with all vertices in $t_{m_1}H$ connecting $x$ and $y$. By Lemma \ref{l3} again, $\alpha$ must lie outside $N_n(H)$. Therefore, the pair $(x,y)$ is either of the form $(h_1t^n, h_2t^n)$ or $(h_1t^{-n}, h_2t^{-n})$.
 
\end{proof}

\begin{prop}
\label{p1}
Let $H$ be a finitely generated group and $G=\langle H,t/tht^{-1}=\phi(h)\rangle$ where $\phi$ in $Aut(H)$. Then, $Div(G,H)\preceq Dist^H_G$. 
\end{prop}
\begin{proof}
Let $T$ be a finite generating set of $H$ and let $S=T\cup\{t\}$. Then, $S$ is a finite generating set of $G$. Suppose that $Div(G,H)=\{\delta^n_{\rho}\}$. We will show that $\delta^n_{\rho}(r)\leq Dist^H_G(nr)$ for all positive integer $r$.

Indeed, let $x$, $y$ be arbitrary points in $\partial N_r(H)$ such that $d_{r,H}(x,y)<\infty$ and $d_S(x,y)\leq nr$. By Lemma \ref{l5}, $x$, $y$ both lie in the same coset $t^mH$ where $\abs{m}=r$. Therefore, there is a path $\gamma$ with all vertices in $t^mH$ connecting $x$ and $y$ and the length of $\gamma$ is at most $Dist^H_G(nr)$. By Lemma \ref{l3} again, the path $\gamma$ must lie outside $N_r(H)$. Therefore, $d_{\rho r,H}(x,y)\leq Dist^H_G(nr)$. Thus, $\delta^n_{\rho}(r)\leq Dist^H_G(nr)$. This implies that $Div(G,H)\preceq Dist^H_G$. 
\end{proof}

\section{Relative divergence of right-angled Artin groups with respect to Bestvina-Brady subgroups and subgroup distortion of Bestvina-Brady subgroups}

From now, we let $\Gamma$ be a finite, connected, simplicial graph with at least 2 vertices. Let $A_{\Gamma}$ be the associated right-angled Artin group and $H_{\Gamma}$ be its Bestvina-Brady subgroup. Let $X_{\Gamma}$ be the associated Salvetti complex and $\tilde{X}_{\Gamma}$ its universal covering. We consider the 1--skeleton of $\tilde{X}_{\Gamma}$ as a Cayley graph of $A_{\Gamma}$ and the vertex set $S$ of $\Gamma$ as a finite generating set of $A_{\Gamma}$. By Theorem \ref{th1}, we can choose the set $T$ of all elements of the form $st^{-1}$ whenever $s$ and $t$ are adjacent vertices as a finite generating set for $H_{\Gamma}$. Let $\Phi$ and $\tilde{f}$ be group homomorphism and continuous map as in Remark \ref{r1}.

\begin{lem}
\label{l2}
Let $M$ be the diameter of $\Gamma$. Let $a$ and $b$ be arbitrary vertices in $S$. For each integer $m$, the length of $a^mb^{-m}$ with respect to $T$ is at most $M\abs{m}$. 
\end{lem}

\begin{proof}
Since the diameter of $\Gamma$ is $M$, we can choose positive integer $n\leq M$ and $n+1$ generators $s_0, s_1, \cdots, s_n$ in $S$ such that the following conditions hold:
\begin{enumerate}
\item $s_0=a$ and $s_n=b$.
\item $s_i$ and $s_{i+1}$ commutes where $i\in \{0, 1, 2, \cdots, n-1\}$.
\end{enumerate} 

Obviously,
\begin{align*}
a^mb^{-m}=s^m_0s^{-m}_n &=(s^m_0s^{-m}_1)(s^m_1s^{-m}_2)(s^m_2s^{-m}_3)\cdots (s^m_{n-2}s^{-m}_{n-1}) (s^m_{n-1}s^{-m}_n)\\ &= (s_0s^{-1}_1)^m(s_1s^{-1}_2)^m(s_2s^{-1}_3)^m\cdots (s_{n-2}s^{-1}_{n-1})^m(s_{n-1}s^{-1}_{n})^m.
\end{align*} 

Also, $s_{i-1}s^{-1}_{i}$ belongs to $T$. Therefore, the length of $a^mb^{-m}$ with respect to $T$ is at most $n\abs{m}$. This implies that the length of $a^mb^{-m}$ with respect to $T$ is at most $M\abs{m}$.
\end{proof}

\begin{prop}
\label{p2}
The subgroup distortion $Dist^{H_{\Gamma}}_{A_{\Gamma}}$ is dominated by a quadratic function. Moreover, $Dist^{H_{\Gamma}}_{A_{\Gamma}}$ is linear when $\Gamma$ is a join. 
\end{prop}

\begin{proof}
We first show that $Dist^{H_{\Gamma}}_{A_{\Gamma}}$ is dominated by a quadratic function. Let $n$ be an arbitrary positive integer and $h$ be an arbitrary element in $H_{\Gamma}$ such that $\abs{h}_S\leq n$. We can write $h=s^{m_1}_1s^{m_2}_2s^{m_3}_3\cdots s^{m_k}_k$ such that:
\begin{enumerate}
\item Each $s_i$ lies in $S$, $\abs{m_i}\geq 1$ and $\abs{m_1}+\abs{m_2}+\abs{m_3}+\cdots+\abs{m_k}\leq n$.
\item $m_1+m_2+m_3+\cdots+m_k=0$
\end{enumerate}

Obviously, we can rewrite $h$ as follows:
\[h=(s^{m_1}_1s^{-m_1}_2)(s^{(m_1+m_2)}_2s^{-(m_1+m_2)}_3)\cdots (s^{(m_1+m_2+\cdots+m_{k-1})}_{k-1}s^{-(m_1+m_2+\cdots+m_{k-1})}_{k}).\]

Let $M$ be the diameter of $\Gamma$. By Lemma \ref{l2}, we have
\begin{align*}
\abs{h}_T &\leq M\abs{m_1}+M\abs{m_1+m_2}+\cdots+M\abs{m_1+m_2+\cdots+m_{k-1}}\\&\leq M\abs{m_1}+M\bigl(\abs{m_1}+\abs{m_2}\bigr)+\cdots+M\bigl(\abs{m_1}+\abs{m_2}+\cdots+\abs{m_{k-1}}\bigr)\\&\leq M(k-1)n\leq Mn^2.
\end{align*}
Therefore, the distortion function $Dist^{H_{\Gamma}}_{A_{\Gamma}}$ is bounded above by $Mn^2$.

We now assume that $\Gamma$ is a join of $\Gamma_1$ and $\Gamma_2$. We need to prove that the distortion $Dist^{H_{\Gamma}}_{A_{\Gamma}}$ is linear. Let $n$ be an arbitrary positive integer and $h$ be an arbitrary element in $H_{\Gamma}$ such that $\abs{h}_S\leq n$. Since $A_{\Gamma}$ is the direct product of $A_{\Gamma_1}$ and $A_{\Gamma_2}$, we can write $h=(a^{m_1}_1a^{m_2}_2\cdots a^{m_k}_k)(b^{n_1}_1b^{n_2}_2\cdots b^{n_\ell}_\ell)$ such that:

\begin{enumerate}
\item Each $a_i$ is a vertex of $\Gamma_1$ and each $b_j$ is a vertex of $\Gamma_2$.
\item $\bigl(\abs{m_1}+\abs{m_2}+\cdots+\abs{m_k}\bigr)+\bigl(\abs{n_1}+\abs{n_2}+\cdots+\abs{n_\ell}\bigr)\leq n$.
\item $(m_1+m_2+\cdots+m_k)+(n_1+n_2+\cdots+n_\ell)=0$
\end{enumerate}
Let $m=m_1+m_2+\cdots+m_k$. Then, $n_1+n_2+\cdots+n_\ell=-m$ and $\abs{m}\leq n$. Let $a$ be a vertex in $\Gamma_1$ and $b$ a vertex in $\Gamma_2$. Since $a$ commutes with each $b_j$, $b$ commutes with each $a_i$ and $a$, $b$ commute, we can rewrite $h$ as follows: 
\begin{align*}
h&=(a^{m_1}_1a^{m_2}_2\cdots a^{m_k}_kb^{-m})(b^ma^{-m})(a^mb^{n_1}_1b^{n_2}_2\cdots b^{n_\ell}_\ell)\\&=(a_1b^{-1})^{m_1}(a_2b^{-1})^{m_2}\cdots (a_kb^{-1})^{m_k}(ba^{-1})^m(ab_1^{-1})^{-n_1}(ab_2^{-1})^{-n_2}\cdots(ab_{\ell}^{-1})^{-n_\ell}.
\end{align*}
Also, $ab_j^{-1}$, $a_ib^{-1}$ and $ba^{-1}$ all belong to $T$. Therefore, 
\begin{align*}
\abs{h}_T\leq \bigl(\abs{m_1}+\abs{m_2}+\cdots+\abs{m_k}\bigr)+\bigl(\abs{n_1}+\abs{n_2}+\cdots+\abs{n_\ell}\bigr)+\abs{m}\leq 2n.
\end{align*}
Therefore, the distortion function $Dist^{H_{\Gamma}}_{A_{\Gamma}}$ is bounded above by $2n$.
\end{proof}

\begin{prop}
\label{p3}
If $\Gamma$ is not a join graph, then the relative divergence $Div(A_{\Gamma}, H_{\Gamma})$ is at least quadratic.
\end{prop}

\begin{proof}
Let $J$ be a maximal join in $\Gamma$ and let $v$ be a vertex not in $J$. Let $g$ in $A_J$ be the product of all vertices in $J$. Let $n=\Phi(g)$ and let $h=gv^{-n}$. Then $h$ is an element in $H_{\Gamma}$. Since $J$ is a maximal join in $\Gamma$ and let $v$ be a vertex not in $J$, then $h$ does not lie in any join subgroup. Also, $h$ is a cyclically reduced element. Therefore, the divergence of the bi-infinite geodesic $\alpha=\cdots hhhhh \cdots$ is at least quadratic by Lemma \ref{l4}. 

Let $t$ be an arbitrary generator in $S$ and $k=\abs{h}_{S}$. We can assume that $\alpha(0)=e$, $\alpha(km)=h^m$, and $\alpha(-km)=h^{-m}$. In order to prove the relative divergence $Div(A_{\Gamma}, H_{\Gamma})$ is at least quadratic, it is sufficient to prove each function $\delta_\rho^n$ dominates the divergence function of $\alpha$ for each $n\geq 2k+2$.

Indeed, let $r$ be an arbitrary positive integer. Let $x=h^{-r}t^r$ and $y=h^{r}t^r$. By the similar argument as in Lemma \ref{l3} and Lemma \ref{l5}, two points $x$ and $y$ both lie in $\partial N_r(H_{\Gamma})$ and $d_{r,H_{\Gamma}}(x,y)<\infty$. Moreover,
\[d_S(x,y)\leq d_S(x,h^{-r})+d_S(h^{-r},h^{r})+d_S(h^{r},y)\leq r+2kr+r\leq (2k+2)r\leq nr.\] 

Let $\gamma$ be an arbitrary path outside $N_{\rho r}(H)$ connecting $x$ and $y$. Obviously, the path $\gamma$ must lie outside the open ball $B\bigl(\alpha(0),\rho r\bigr)$. It is obvious that we can connect $x$ and $h^{-r}$ by a path $\gamma_1$ of length $r$ which lies outside $B\bigl(\alpha(0),\rho r\bigr)$. Similarly, we can connect $y$ and $h^{r}$ by a path $\gamma_2$ of length $r$ which lies outside $B\bigl(\alpha(0), \rho r\bigr)$. Let $\gamma_3$ be the subsegment of $\alpha$ connecting $\alpha(-\rho r)$ and $h^{-r}$. Let $\gamma_4$ be the subsegment of $\alpha$ connecting $\alpha(\rho r)$ and $h^{r}$. It is not hard to see the length of $\gamma_3$ and $\gamma_4$ are both $(k-\rho)r$.

Let $\bar{\gamma}=\gamma_3\cup\gamma_1\cup \gamma\cup \gamma_2\cup\gamma_4$. Then, $\bar{\gamma}$ is a path that lies outside $B\bigl(\alpha(0), \rho r\bigr)$ connecting $\alpha(-\rho r)$ and $\alpha(\rho r)$. Therefore, the length of $\bar{\gamma}$ is at least $Div_{\alpha}(\rho r)$.
Also,
\[\ell(\bar{\gamma})=\ell(\gamma_3)+\ell(\gamma_1)+\ell(\gamma)+\ell(\gamma_2)+\ell(\gamma_5)=\ell(\gamma)+2(k-\rho+1)r.\]
Thus, 
\[\ell(\gamma)\geq Div_{\alpha}(\rho r)-2(k-\rho+1)r.\] 
This implies that \[d_{\rho r,H_{\Gamma}}(x,y)\geq Div_{\alpha}(\rho r)-2(k-\rho+1)r.\]

Therefore, \[\delta_\rho^n(r)\geq Div_{\alpha}(\rho r)-2(k-\rho+1)r.\]
Thus, the relative divergence $Div(A_{\Gamma}, H_{\Gamma})$ is at least quadratic.
\end{proof}

The following theorem is deduced from Proposition \ref{p1}, Proposition \ref{p2}, and Proposition \ref{p3}. 
\begin{thm}
Let $\Gamma$ be a connected, finite, simplicial graph with at least 2 vertices. Let $A_{\Gamma}$ be the associated right-angled Artin group and $H_{\Gamma}$ the Bestvina-Brady subgroup. Then the relative divergence $Div(A_{\Gamma}, H_{\Gamma})$ and the subgroup distortion $Dist^{H_{\Gamma}}_{A_{\Gamma}}$ are both linear if $\Gamma$ is a join graph. Otherwise, the relative divergence $Div(A_{\Gamma}, H_{\Gamma})$ and the subgroup distortion $Dist^{H_{\Gamma}}_{A_{\Gamma}}$ are both quadratic. 
\end{thm}

\begin{cor}
For each integer $n\geq 3$, there is a right-angled Artin group containing a free subgroup of rank $n$ whose inclusion is not a quasi-isometric embedding.
\end{cor}

\begin{proof}
For each positive integer $n\geq 3$, let $\Gamma$ be a tree with $n$ edges such that $\Gamma$ is not a join graph. By the above theorem, the distortion of $H_\Gamma$ in the right-angled Artin group $A_\Gamma$ is quadratic. Also, $H_\Gamma$ is the free group of rank $n$ by Theorem \ref{th1}.
\end{proof}

\section{Connection to horizontal surface subgroups}

\begin{defn}
A \emph{graph manifold} is a compact, irreducible, connected orientable $3$--manifold $M$ that can be decomposed along $\mathcal{T}$ into finitely many Seifert manifolds, where $\mathcal{T}$ is the canonical decomposition tori of Johannson and of Jaco-Shalen. We call the collection $\mathcal{T}$ is JSJ-decomposition in $M$, and each element in $\mathcal{T}$ is JSJ-torus.
\end{defn}

%Define an associated graph $\Gamma(M)$ as below: Each Seifert component $P_v$ determines a vertex $v$ and each JSJ-torus $T_e$ determines an edge $e$, vertices $v$ and $v'$ are connected by an edge $e$ if and only if $T_e$ is shared by $P_v$ and $P_{v'}$. If we put an orientation on $e$ for each $e$ in $\Gamma(M)$, then $e$ determines a homeomorphism $g_e\!: T_e\to \bar{T}_e$, where $T_e$ and $\bar{T}_e$ are tori in the boundaries of the vertex manifolds, corresponding to the beginning and the end of $e$ respectively.

\begin{defn}
If $M$ is a Seifert manifold, a properly immersed surface $g\!: S \looparrowright M$ is \emph{horizontal}, if $g(S)$ is transverse to
the Seifert fibers everywhere. In the case $M$ is a graph manifold, a properly immersed surface $g\!: S \looparrowright M$ \emph{horizontal} if $g(S)\cap P_v$ is horizontal for every Seifert component $P_v$.
\end{defn}

\begin{figure}%{.5\textwidth}
\centering
\labellist
\small\hair 2pt

\pinlabel $u_1$ at 300 415
\pinlabel $u_2$ at 400 480
\pinlabel $u_3$ at 520 440
\pinlabel $u_4$ at 425 370
\pinlabel $u_5$ at 325 285
\pinlabel $u_6$ at 545 300
\pinlabel $u_7$ at 425 230
\pinlabel $b_v$ at 425 570

\endlabellist
%\scalebox{.2}{\includegraphics{Graph}}
\includegraphics[scale=0.35]{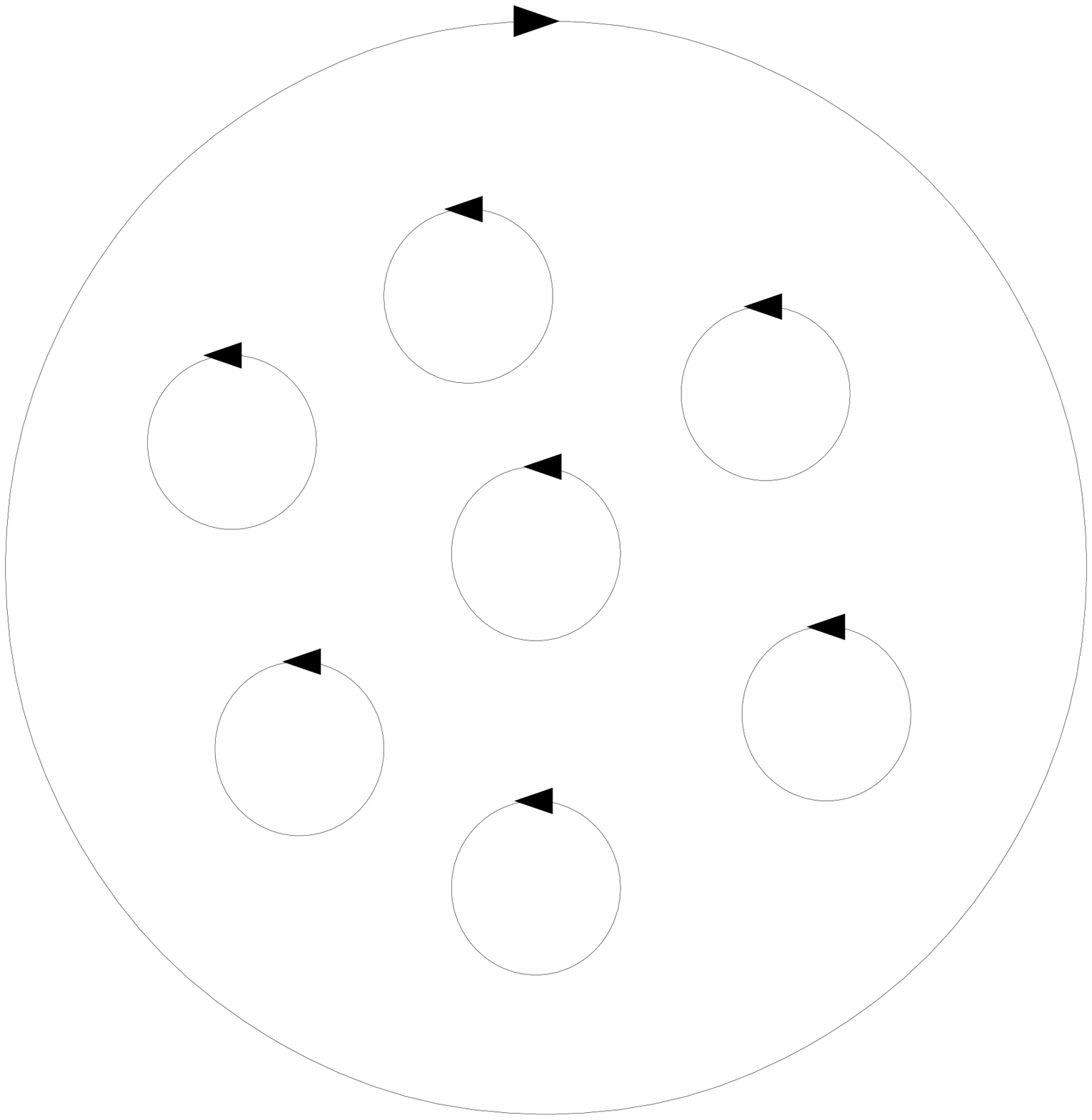}
\caption{A punctured disk $\Sigma_v$ when the degree of $v$ in $\Gamma$ is 7.}
\label{i1}
%\end{subfigure}
%\caption{}
%\label{a}
\end{figure}

\begin{thm}
If $\Gamma$ is a finite tree with at least 3 vertices, then the associated right-angled Artin group $A_\Gamma$ is a fundamental group of a graph manifold and the Bestvina-Brady subgroup $H_\Gamma$ is a horizontal surface subgroup.
\end{thm}

\begin{proof}
First, we construct the graph manifold $M$ whose fundamental group is $A_{\Gamma}$ as follows:

Let $v$ be a vertex of $\Gamma$ of degree $k\geq 2$. Let $u_1, u_2,\cdots,u_k$ be all elements in $\ell k(v)$. % Obviously $A_{st(v)}=F_k \times \Z$, where $k$ is the number of points in the link $lk(v)$ of $v$, $F_k$ is the free group generated by element in $lk(v)$ and the factor $\Z$ generated by $v$. 
Let $\Sigma_v$ be a punctured disk with $k$ inside holes in which their boundaries are labeled by elements in $\ell k(v)$. We also label the outside boundary component of $\Sigma_v$ by $b_v$ (see Figure \ref{i1}). Obviously, $\pi_1(\Sigma_v)$ is the free group generated by $u_1, u_2,\cdots u_k$.

Let $P_v=\Sigma_v\times S^1_v$, here we label the circle factor in $P_v$ by $v$. Obviously, each $P_v$ is a Seifert manifold. Moreover, for each $u_i$ in $\ell k(v)$, the Seifert manifold $P_v$ contains torus $S^1_{u_i}\times S^1_v$ as a component of its boundary. 

We construct the graph manifold by gluing pair of Seifert manifolds $(P_{v_1},P_{v_2})$ along their tori $S^1_{v_1}\times S^1_{v_2}$ whenever $v_1$ and $v_2$ are adjacent vertices in $\Gamma$. We observe that the pair of such regions are glued together by switching fiber and base directions. It is not hard to see that the fundamental group of $M$ is the right-angled Artin group $A_\Gamma$.

We now construct the horizontal surface $S$ in $M$ with the Bestvina-Brady subgroup $H_\Gamma$ as its fundamental group. We first construct the horizontal surface $S_v$ on each Seifert piece $P_v=\Sigma_v\times S^1_v$, where $v$ is a vertex of $\Gamma$ of degree $k\geq 2$. 

We remind the reader that $\Sigma_v$ is a punctured disk with $k$ inside holes in which their boundaries are labeled by all elements $u_1, u_2, \cdots,u_k$ in $\ell k(v)$. We also label the outside boundary component of $\Sigma_v$ by $b_v$ (see Figure \ref{i1}). We label the circle factor in $P_v$ by $v$. %For each $u_i$ in $\ell k(v)$, we construct a simple curve $c_{u_i,v}=[u_i]-[v]$ in the boundary torus $S^1_{u_i}\times S^1_v$. 
%In the boundary torus $S^1_{b_v}\times S^1_v$, we construct a simple curve $c_{b_v,v}=[b_v]+k[v]$.

%We claim that the family of curves $c_{u_1,v}, c_{u_2,v}, \cdots, c_{u_k,v}, c_{b_v,v}$ bounds a horizontal surface $S_v$ of genus 0. Indeed, 
Let $S_v$ be a copy of the punctured disk $\Sigma_v$. However, we relabel all inside circles by $c_1, c_2,\cdots, c_k$ and the outside circle by $c_v$. We will construct a map $(g,h)\!: S_v\to \Sigma_v\times S^1_v$ as follows:
\begin{enumerate}
\item The map $g$ is the identity map that maps each $c_i$ to $u_i$ and $c_v$ to $b_v$.
\item The map $h$ has degree $-1$ on boundary component $c_i$ and degree $k$ on $c_v$.
\end{enumerate}

We now construct the map $h$ with the above properties. We observe that the fundamental group of $S_v$ is generated by $c_1, c_2,\cdots, c_k$, and $c_v$ with a unique relator $c_1 c_2 c_3\cdots c_k c_v =e$. Here we abused notation for the presentation of $\pi_1(S_v)$. By that presentation of $\pi_1(S_v)$, we can see that there is a group homomorphism $\phi$ from $\pi_1(S_v)$ to $\Z$ that maps each $c_i$ to $-1$ and $c_v$ to $k$. By \cite[Proposition 1B.9]{MR1867354}, the group homomorphism $\phi$ is induced by a map $h$ from $S_v$ to $S^1_v$. Therefore, we constructed a desired map $h$.

%We note that the map $h$ does exist because the sum of degrees on all boundary components of $S_v$ is $0$, the fact that $[S_v, S^1_v] = H^1(S_v;\Z)$ and the exact cohomology sequence 
%\[H^1(S_v;\Z) \to H^1(\partial S_v;\Z) \to H^2(S_v, \partial S_v;\Z) = \Z.\]

Finally, we identify the surface $S_v$ with its image via the map $(g,h)$. By construction, $\pi_1(S_v)$ is the subgroup of $\pi_1(P_v)$ generated by elements $u_1v^{-1}, u_2v^{-1}, \cdots, u_kv^{-1}$. We observe that if we glue pair of Seifert manifolds $(P_{v_1},P_{v_2})$ along their tori $S^1_{v_1}\times S^1_{v_2}$, pair of horizontal surfaces $(S_{v_1},S_{v_2})$ will be matched up along their boundaries in $S^1_{v_1}\times S^1_{v_2}$. Therefore, we constructed a horizontal surface $S$ in $M$. By Vankampen theorem, the fundamental group of $S$ is generated by all elements of the form $st^{-1}$ whenever $s$ and $t$ are adjacent vertices in $\Gamma$. In other words, $\pi_1(S)$ is the Bestvina-Brady subgroup by Theorem \ref{th1}.
\end{proof}

%%%%%%%%%%%%%%%%%%%%%%%%%%%%%%%%%%%%%%%%%%%%%%%%%%%%%%%%%%%
%%                BIBLIOGRAPHY
%%%%%%%%%%%%%%%%%%%%%%%%%%%%%%%%%%%%%%%%%%%%%%%%%%%%%%%%%%%
\bibliographystyle{alpha}
\bibliography{Tran}
\end{document}